\documentstyle[amssymb,12pt]{amsart}

\ifx\mathrm\undefined\let\mathrm\bf\fi
\ifx\mathbf\undefined\let\mathbf\bf\fi

 \def
\Sum{\sum\limits}

\let\dl\delta \let\Dl\Delta
\let\epe\epsilon \let\eps\varepsilon \let\epsilon\eps

\let\la\lambda 
 
 \let\phi\varphi

\newcommand{\half}{\frac12}
\newcommand{\Z}{{\Bbb Z}}  
\newcommand{\N}{{\Bbb N}}

\newcommand{\C}{{\Bbb C}}
   
\newcommand{\Ref}[1]{{$($\ref{#1}$)$}}
\newcommand{\bean}{\begin{eqnarray}}
\newcommand{\eean}{\end{eqnarray}}
\newcommand{\be}{\begin{displaymath}}
\newcommand{\ee}{\end{displaymath}}
\newcommand{\bea}{\begin{eqnarray*}}   
\newcommand{\eea}{\end{eqnarray*}}
\newcommand{\g}{{{\frak g}\,}}

\newcommand{\h}{{{\frak h\,}}}

\newcommand{\T}{\otimes}

\newcommand{\vs}{\vspace{1.5\baselineskip}}

\newenvironment{proof}{\noindent{\it Proof\/}:\rm}{$\;\Box$
\par\vs}

\newtheorem
{thm}{Theorem}

\newtheorem
{lemma}[thm]{Lemma}

\newcommand{\th}{\theta}
\newcommand{\End}{{\operatorname{End}}}

\begin{document}

\title{Quantization of the space of conformal blocks.} 
\author{E. Mukhin and A. Varchenko}
\maketitle
\vskip-.5\baselineskip
\centerline{\it Department of Mathematics,
University of North Carolina at Chapel Hill,}
\centerline{\it Chapel Hill, NC 27599-3250, USA}
\centerline{{\it E-mail addresses:} {\rm mukhin@@math.unc.edu,
av@@math.unc.edu}}
\bigskip
\medskip
\centerline{October, 1997}
\bigskip
\medskip

\begin{abstract}
We consider the discrete Knizhnik-Zamolodchikov connection (qKZ) associated to $gl(N)$, defined in terms of rational R-matrices. We prove that under certain resonance conditions, the qKZ connection has a non-trivial invariant subbundle which we call the subbundle of quantized conformal blocks. The subbundle is given explicitly by algebraic equations in terms of the Yangian $Y(gl(N))$ action. The subbundle is a deformation of the subbundle of conformal blocks in CFT. The proof is based on an identity in the algebra with two generators $x,y$ and defining relation $xy=yx+yy$. 
\end{abstract}

\section{Introduction}

Conformal Field Theory (CFT) associates a finite dimensional vector space, called the space of conformal blocks, to each Riemann surface with marked points and certain additional data (local coordinates, representations). The vector spaces corresponding to different complex structures or different positions of the points are locally (projectively) identified by a projectively flat connection.

A Wess-Zumino-Witten (WZW) model is labeled by a simple Lie algebra $\g$ and a positive integer $c$ called level. The space of conformal blocks is defined in terms of the representation theory of the affine Kac-Moody algebra $\widehat{L\g}$, which is a central extension of the loop algebra $L\g=\g\T\C((t))$, see \cite{K}. For each irreducible highest weight $\g$-module $L$, we have a canonically defined corresponding irreducible highest weight $\widehat{L\g}$ module of level $c$, denoted $\widehat{L}$, see \cite{K}. Then, the space of conformal blocks associated to a Riemann surface $\Sigma$, $n$ distinct points with a choice of local holomorphic parameters around them and $n$ irreducible finite dimensional $\g$ modules $L_1,\dots, L_n$ (such that $\widehat{L}_i$ are integrable) is the space $H^0(\g(\Sigma),(\widehat{L}_1\T\ldots\T \widehat{L}_n)^*)$ of linear forms on $\widehat{L}_1\T\ldots\T \widehat{L}_n$ invariant under the action of the Lie algebra $\g(\Sigma)$ of meromorphic $\g$-valued functions on $\Sigma$, holomorphic outside of the marked points. The action of $\g(\Sigma)$ is defined through the Laurent expansion at the poles. Varying the data makes the spaces of conformal blocks into a holomorphic vector bundle with a projectively flat connection given by the Sugawara-Segal construction.

An explicit description is known on the Riemann sphere ${\Bbb P}^1$. Namely, the space of conformal blocks on ${\Bbb P}^1$ with $n+1$ distinct points $z_1,\dots,z_n\in\C\subset{\Bbb P}^1$, $z_{n+1}=\infty$, associated to $n+1$ irreducible $\g$ modules $L_1,\dots,L_{n+1}$ with highest weights $\la_1,\dots,\la_{n+1}$, is identified with a subspace of $(L_1\T\ldots\T L_n)^{sing}_\la$, where  $(L_1\T\ldots\T L_n)^{sing}_\la\subset L_1\T\ldots\T L_n$ is the weight subspace of singular vectors of total weight $\la$, $\la=-w\la_{n+1}$ and $w$ is the longest element of the Weyl group. More precisely the space of conformal blocks can be described as follows. Let $\theta$ be the highest root and let the scalar product be normalized by $(\theta,\theta)=2$.

If the resonance condition, 
\bean\label{res'}
(\theta,\la)-c+k-1=0,
\eean
holds for some $k\in\N$, then the space of conformal blocks is identified with subspace
\be
W_{\la_{n+1}}(z)=\{m\in (L_1\T\ldots\T L_n)^{sing}_\la\,|\, (E(z))^km=0\},
\ee
otherwise the space of conformal blocks is identified with the weight space of singular vectors,
$W_{\la_{n+1}}(z)=(L_1\T\ldots\T L_n)^{sing}_\la$.
Here
\be
E(z)=\sum_{i=0}^n z_ie_\th^{(j)},
\ee
and $e_\th^{(j)}$ denotes the element $e_\th\in \g$ acting on the $j$-th factor, see \cite{KT}, \cite{FSV1}.

The subspaces $W_{\la_{n+1}}(z)$ associated to different sets of marked points $z_1,\dots,z_n\in\C$ form a holomorphic subbundle of the trivial vector bundle over the configuration space $\C^{[n]}=\{z\in\C^n\,|\, z_i\neq z_j,\, i,j=1,\dots,n,\, i<j\}$ with fiber $L_1\T\ldots\T L_n$. There is a flat connection on the bundle $\C^{[n]}\times (L_1\T\ldots\T L_n)$ preserving the subbundle of conformal blocks. Its horizontal sections $\Psi(z)$ obey the Knizhnik-Zamolodchikov equation, 
\bean\label{KZ}
\partial_{z_i}\Psi(z)=\frac{1}{\kappa}\sum_{j\neq i}\frac{\Omega_{ij}}{z_i-z_j}\Psi(z),\qquad i=1,\dots,n,
\eean
$\kappa=c+h^\vee$, where $\Omega_{ij}$ is the Casimir operator acting on the $i$-th and $j$-th factors and $h^\vee$ is the dual Coxeter number of $\g$.

Let $\g=gl(N)$. The quantized Knizhnik-Zamolodchikov (qKZ) equation is a holonomic system of difference equations for a function $\Psi(z)$ with values in $L_1\T\ldots\T L_n$,
\be
\Psi(z_1,\dots,z_i+p,\dots,z_n)=R_{m,m-1}(z_m-z_{m-1}+p)\ldots R_{m,1}(z_m-z_1+p)\times
\ee
\be
\times
R_{m,n}(z_m-z_n)\ldots R_{m,m+1}(z_m-z_{m+1})\Psi(z),
\ee
$i=1,\dots,n$, where $R_{i,j}(x)$ is the rational $R$-matrix $R_{L_{\la_i}L_{\la_j}}(x)$ acting in $i$-th and $j$-th factors and $p\in\C$ is a parameter, see \cite{FR}.

The qKZ equation defines the discrete Knizhnik-Zamolodchikov (qKZ) connection on the trivial vector bundle over $\C^n$ with fiber $L_1\T\ldots\T L_n$.
 
Consider the quasiclassical limit of the qKZ. Namely, set $y_i=z_i/h$, for some $h\in\C$ and let $h\to 0$. In this limit the qKZ equation turns into a system of differential equations
\bean\label{KZ'}
p\partial_{y_i}\tilde{\Psi}(y)=-\sum_{j\neq i}\frac{\tilde{\Omega}_{ij}}{y_i-y_j}\tilde{\Psi}(y),\qquad i=1,\dots,n,
\eean
where $\tilde{\Omega}_{ij}=\Omega_{ij}+A_{ij}$ and $A_{ij}\in\C$ is a constant defined by 
\be
\Omega v_i\T v_j=-A_{ij}v_i\T v_j,
\ee
$v_i,v_j$ are highest weight vectors generating $L_i,L_j$, see Section 7 in \cite{TV1} and Section 12.5 in \cite{CP}.

Notice that if $p=-\kappa$, then $\Psi(z)$ is a solution of the KZ equation \Ref{KZ} if and only if the function
\be
\tilde{\Psi}(y)=\prod_{i<j}(y_i-y_j)^{A_{ij}/\kappa}\,\Psi(y)
\ee
is a solution of the equation \Ref{KZ'}.

In this paper we suggest a quantization of the space of conformal blocks. Namely, under the resonance condition, 
\bean\label{res}
(\theta,\la)+p+N+k-1=0,
\eean
where $k\in\N$, we introduce the space of quantized conformal blocks, $C_{\la_{n+1}}(z)$, by
\be
C_{\la_{n+1}}(z)=\{m\in (L_1\T\ldots\T L_n)^{sing}_\la\,|\, (e(z))^km=0\},
\ee
where
\be
e(z)m=\sum_{j=1}^n \left(z_j-e_{NN}^{(j)}+\sum_{s=j+1}^n2h_\th^{(s)}\right)e_\th^{(j)}m+\sum_{j=2}^{N-1}\sum_{r,s=1\atop r<s}^Ne_{Nj}^{(r)}e_{j1}^{(s)}m,
\ee
and $e_{ij}^{(s)}$ denotes the element $e_{ij}\in gl(N)$ acting on the $s$-th factor of $L_1\T\ldots\T L_n$. 

The operator $e(z)$ can be described in terms of the action of the Yangian $Y(gl(N))$ in the tensor product of evaluation modules $L_1(z_1)\T\ldots\T L_n(z_n)$, $e(z)=T_{N1}^{(2)}-T_{NN}^{(1)}T_{N1}^{(1)}$, see \Ref{e(z)-Yangian}.

In the quasiclassical limit the resonance condition \Ref{res} coincides with the resonance condition \Ref{res'} and the operator $e(z)$ tends to the operator $E(z)$.

In this paper we show that the space of quantized conformal blocks is invariant with respect to the quantized Knizhnik-Zamolodchikov connection. The proof is based on an identity in the algebra with two generators $x,y$ and defining relation $xy=yx+yy$.

In a recent paper \cite{EF}, B.Enriquez and G.Felder gave a construction of the space of quantized conformal blocks as a space of suitable coinvariants of an action of quantum doubles of Yangians. We expect that the Enriquez-Felder construction applied to our situation will identify the space of coinvariants with the space of conformal blocks introduced in this paper.

The simplest qKZ equation is the qKZ equation associated with $sl(2)$. In \cite{TV1}, \cite{MV} solutions of the $sl(2)$ qKZ equation were constructed in terms of multidimensional hypergeometric functions. In the next paper we will show that all hypergeometric solutions of the $sl(2)$ qKZ equation automatically belong to the space of quantized conformal blocks. This result is analogues to the fact that the values of all hypergeometric solutions of the KZ differential equation automatically belong to the space of conformal blocks of CFT, see [FSV1-3].

\section{The qKZ connection}
\subsection{The Lie algebra $gl(N)$.}
Let $N$ be a natural number. The Lie algebra $gl(N)$ is spanned over $\C$ by elements $e_{ij}$, $i,j=1,\dots,N$, with commutators given by
\be
[e_{ij},e_{kl}]=\dl_{jk}e_{il}-\dl_{il}e_{kj}, \qquad i,j,k,l=1,\dots,N,
\ee
where $\dl_{ij}$ is the Kronecker symbol.
 
Let $\h$ be Cartan subalgebra spanned by $h_i=\half e_{i,i}$, $i=1,\dots,N$. Let $\th$ be the highest root and let $h_\th=h_1-h_N$, $e_\th=e_{1,N}$.

Let $\la=(\la_{(1)},\dots,\la_{(N)})\in\C^N$. Let $V$ be the Verma module over $gl(N)$ with highest weight $\la$, $V$, i.e $V$ is generated by a highest vector $v$ such that $h_iv=\la_{(i)}v$ and $e_{ij}v=0$ for $i,j=1,\dots,N$, $i<j$.

Let $S\in V$ be the maximal proper submodule. Then $L=V/S$ is the irreducible $gl(N)$ module with highest weight $\la$.

For a $gl(N)$ module $M$ with highest weight $\la=(\la_{(1)},\dots,\la_{(N)})\in\C^N$, let $M^{sing}\in M$ be the subspace of singular vectors, i.e. the subspace of vectors annihilated by $e_{i,i+1}$, for all $i=1,\dots,N-1$,
\be
M^{sing}=\{m\in M\,|\, e_{i,i+1}m=0,\, i=1,\dots,N-1\}.
\ee
Let also 
\be
(M)_l=\{m\in M\,|\,h_\th m=(\la_{(1)}-\la_{(N)}-l)\,m\},
\ee
and $(M)_l^{sing}=(M)_l\bigcap M^{sing}$.

\subsection{The Hopf algebra $Y(gl(N))$.}
The Yangian $Y(gl(N))$ is an associative algebra with an infinite set of generators $T_{i,j}^{(s)}$,  $i,j=1,\dots,N$, $s=0,1,\dots$, subject to the following relations:
\be
[T_{ij}^{(r)},T_{kl}^{(s+1)}]-[T_{ij}^{(r+1)},T_{kl}^{(s)}]=T_{kj}^{(r)}T_{il}^{(s)}-
T_{kj}^{(s)}T_{il}^{(r)},\qquad T_{ij}^{(0)}=\dl_{ij},
\ee
$i,j,k,l=1,\dots,N$; $r,s,=1,2,\dots$ .

The comultiplication $\Dl\, :\,Y(gl(N))\to Y(gl(N))\T Y(gl(N))$ is given by 
\be
\Dl\,:\,T_{ij}^{(s)}\mapsto \sum_{k=1}^N\sum_{r=0}^s\, T_{ik}^{(r)}\T T_{kj}^{(s-r)}.
\ee 

For each $x\in\C$, there is an automorpfism $\rho_x\,:\, Y(gl(N))\to Y(gl(N))$ given by
\be
\rho_x\,:\,T_{ij}^{(s)}\mapsto\sum_{r=1}^s\,{s-1 \choose r-1}x^{s-r}T_{ij}^{(r)}.
\ee

There is also an \emph{evaluation morphism} $\epe$ to the universal enveloping algebra of $gl(N)$,
$\epe\,:\, Y(gl(N))\to U(gl(N))$, given by
\be
\epe\,:\,T_{ij}^s\mapsto \dl_{1s}e_{ji},
\ee
for $s=1,2,\dots$ .

Introduce the generating series $T_{ij}(u)=\Sum_{s=0}^\infty T_{ij}^{(s)}u^{-s}$. In terms of these series the relations in the Yangian take the form
\be
R(x-y)T_{(1)}(x)T_{(2)}(y)=T_{(2)}(y)T_{(1)}(x)R(x-y),
\ee
where $R(x)=(x\operatorname{Id}+P)\in\End(\C^N\T\C^N)$, $P\in\End(\C^N\T\C^N)$ is the operator of permutation of the two factors, 
$T_{(1)}(x)=1\T T(x)$, $T_{(2)}(x)=T(x)\T 1$.

In terms of the generating series the comultiplication $\Dl$, the automorphisms $\rho_x$ and the evaluation morphism $\epe$ take the form
\be
\Dl\,:\, T_{ij}(u)\mapsto \sum_{k=1}^N T_{ik}(u)\T T_{kj}(u),
\ee
\be
\rho_x\,:\,T_{ij}(u)\mapsto T_{ij}(u-x),
\ee
\be
\epe\,:\, T_{ij}(u)\mapsto \dl_{ij}+e_{ji}u^{-1},
\ee
$i,j=1,\dots,N$. For more detail on the Yangian see \cite{CP},\cite{KR}.

For any $gl(N)$ module $M$ and $x\in\C$, let $M(x)$ be the $Y(gl(N))$ module obtained from the module $M$ via the homomorphism $\epe\circ\rho(x)$. The module $M(x)$ is called the \emph{evaluation module}. The action of $Y(gl(N))$ in the evaluation module $M(x)$ is given by
\be
T_{ij}^{(s)}m=x^{s-1}e_{ji}m,
\ee 
for all $m\in M$, $i,j=1,\dots,N$, $s=1,2,\dots$ .

Let $L_1,L_2$ be irreducible $gl(N)$ modules. For generic complex numbers $x,y$, the $Y(gl(N))$ modules $L_1(x)\T L_2(y)$ and $L_2(y)\T L_1(x)$ are irreducible and isomorphic. There is a unique intertwiner of the form $PR_{L_1L_2}(x-y)$ mapping $v_1\T v_2$ to $v_2\T v_1$, where $P$ is the operator of permutation of the two factors and $v_i$ are highest weight vectors generating $L_i$, $i=1,2$. The operator $R_{L_1L_2}(x)\in\End(L_1\T L_2)$ is called the \emph{rational $R$-matrix}, see \cite{CP}, \cite{D}.

Let $L_1,L_1,L_3$ be irreducible $gl(N)$ modules. The rational $R$-matrix satisfies the \emph{Yang-Baxter equation} in $\End(L_1\T L_2\T L_3)$ :
\bean\label{YBE}
R_{L_1L_2}(x-y)R_{L_1L_3}(x)R_{L_2L_3}(y)=R_{L_2L_3}(y)R_{L_1L_3}(x)R_{L_1L_2}(x-y),
\eean
the symmetry relation:
\bean\label{sym}
PR_{L_1L_2}(x)=R_{L_2L_1}(x)P\in\End(L_1\T L_2),
\eean
and the inversion relation:
\bean\label{inverse}
R_{L_1L_2}(x)R_{L_1L_2}(-x)=1\in\End(L_1\T L_2).
\eean

For all $x\in\C$, the rational $R$-matrix $R_{L_1L_2}(x)$ commutes with the action of $gl(N)$ in $L_1\T L_2$ and, in particular, preserves the weight decomposition.

\subsection{The qKZ connection.}\label{qKZ}
Fix a non-zero complex number $p$. Let $L_i$, $i=1,\dots,n$, be irreducible $gl(N)$ modules. For $m=1,\dots,n$, $z=(z_1,\dots,z_n)\in\C^n$, define the \emph{Knizhnik-Zamolodchikov operators} $K_m(z)\in\End(L_1\T\ldots\T L_n)$ by the formula
\bean\label{KZ operators}
\lefteqn{
K_m(z)=
R_{L_mL_{m-1}}(z_m-z_{m-1}+p)\ldots R_{L_mL_1}(z_m-z_1+p)\times}
\\&&
\times
R_{L_mL_n}(z_m-z_n)\ldots R_{L_mL_{m+1}}(z_m-z_{m+1}).
\eean

The KZ operators commute with the action of $gl(N)$, $K_m(z)e_{ij}=e_{ij}K_m(z)$ for all $z\in\C^n$, $i,j=1,\dots,N$ and $m=1,\dots,n$. In particular, the KZ operators preserve the subspaces $(L_1\T\ldots\T L_n)_l$ and $(L_1\T\ldots\T L_n)_l^{sing}$ for all $l\in\Z_{\ge 0}$.

\begin{lemma}\label{compatibility}
The qKZ operators $K_m(z)$ and the rational $R$-matrices $R_{L_iL_{i+1}}(x)$ satisfy the following relations
\be
K_{i+1}(z_1,\dots,z_{i+1},z_i,\dots,z_n)P_{L_iL_{i+1}}R_{L_iL_{i+1}}(z_i-z_{i+1})=P_{L_iL_{i+1}}R_{L_iL_{i+1}}(z_i+p-z_{i+1})K_i(z),
\ee
\be
K_k(z_1,\dots,z_j+p,\dots,z_n)K_j(z_1,\dots,z_n)=K_j(z_1,\dots,z_k+p,\dots,z_n)K_k(z_1,\dots,z_n),
\ee
for all, $k,j=1,\dots,n$ and $i=1,\dots,n-1$.
\end{lemma}

Lemma~\ref{compatibility} follows from properties \Ref{YBE}-\Ref{inverse} of the $R$-matrix.

The operators $K_m(z)$, $m=1,\dots,n$, define a discrete flat connection on the trivial vector bundle over $\C^n$ with fiber $L_1\T\ldots\T L_n$. This connection is called the \emph{quantized Knizhnik-Zamolodchikov connection}. For all $l\in\Z_{\ge 0}$, the subspaces $(L_1\T\ldots\T L_n)_l$ and $(L_1\T\ldots\T L_n)_l^{sing}$ are invariant under the qKZ connection.

A subspace $(L_1\T\ldots\T L_n)^{sing}_l$ is called a \emph{resonance subspace} if for some $k\in\N$, $k\le l$, 
\bean\label{resonance}
2h_\th+(p+N+k-1)\operatorname{Id}=0
\eean
in $(L_1\T\ldots\T L_n)^{sing}_l$.

Set 
\bean\label{e(z)-Yangian}
e(z)=T_{N1}^{(2)}-T_{NN}^{(1)}T_{N1}^{(1)}\in Y(gl(N)).
\eean

For each $m\in L_1\T\ldots\T L_n$, we have
\bean\label{e(z)}
e(z)m=\sum_{j=1}^n \left(z_j-e_{NN}^{(j)}+\sum_{s=j+1}^n2h_\th^{(s)}\right)e_\th^{(j)}m+\sum_{j=2}^{N-1}\sum_{r,s=1\atop r<s}^Ne_{Nj}^{(r)}e_{j1}^{(s)}m,
\eean
where $e_{ij}^{(s)}$ denotes $e_{ij}$ acting on the $s$-th factor of $L_1\T\ldots\T L_n$. 

Introduce the \emph{subspace of quantized conformal blocks} $C(z)\subseteq (L_1\T\ldots\T L_n)^{sing}_l$. For a resonance subspace $(L_1\T\ldots L_n)^{sing}_l$, let $C(z)$ be the kernel of the operator $(e(z))^k$ acting in $(L_1\T\ldots\T L_n)^{sing}_l$.
\be
C(z)=\{m\in (L_1\T\ldots\T L_n)^{sing}_l\,|\,(e(z))^km=0\}.
\ee
For a non-resonance subspace $(L_1\T\ldots\T L_n)_l^{sing}$, let $C(z)=(L_1\T\ldots\T L_n)^{sing}_l$. 

\begin{thm}\label{subbundle}
Let $L_{\la_i}$, $i=1,\dots,n$, be irreducible $gl(N)$ modules.  Then the space of conformal blocks $C(z)$ is invariant with respect to the qKZ connection,  
\be
K_i(z)C(z)=C(z_1,\dots,z_i+p,\dots,z_n),
\ee
as well as with respect to permutations of variables, 
\be
PR_{L_iL_{i+1}}(z_i-z_{i+1})C(z)=C(z_1,\dots,z_{i+1},z_i,\dots,z_n).
\ee
\end{thm}

Theorem~\ref{subbundle} is proved in Section~\ref{proof}.

{\bf Remark.} It follows from the proof that the supspace of conformal blocks $C(z)$ in a resonance subspace  $(L_1\T\ldots\T L_n)^{sing}_l$ can be also defined as the kernel of the operator $(T_{N1}^{(2)})^k$ acting in $(L_1\T\ldots\T L_n)^{sing}_l$.

\subsection{The Jordan plane.}
The proof of Theorem~\ref{subbundle} is based on an identity which holds in the algebra called the Jordan plane.

The associative algebra with generators $x,y$ subject to the relation $xy=yx+yy$ is called the \emph{Jordan plane} and is denoted $A^J$, $A^J=\C\langle x,y\rangle/(xy-yx-yy)$, see \cite{DMMZh}.

Let $q\in\C$, $q\neq 0$. The associative algebra with generators $x,y$ subject to the relation $xy=qyx$ is called the \emph{quantum plane with parameter q} and is denoted $A_q$, $A_q=\C\langle x,y\rangle/(xy=qyx)$. The quantum plane with parameter 1 is isomorphic to the ring of polynomials in commuting variables $x,y$, $A_1\simeq\C[x,y]$.

An associative algebra $A$ is called \emph{quadratic} if it has the form $A=\C\langle x_1,\dots,x_s\rangle/(I)$, and the ideal $I$ is generated by homogenious polynomials in $x_1,\dots,x_s$ of degree two. We have $A=\bigoplus\limits_{r=0}^\infty (A)_r$, where  $(A)_r$ is the $r$-th homogenious component spanned over $\C$ by homogenious polynomials in $x_1,\dots,x_r$ of degree $r$.

A quadratic algebra $A$ with two generators has a \emph{polynomial growth} if $\operatorname{dim}(A)_r=\operatorname{dim}\C[x,y]_r=r+1$.

\begin{thm}\label{Jordan} Any quadratic algebra in two generators with a polynomial growth is isomorphic either to the Jordan plane $A^J$ or to the quantum plane $A_q$ with some parameter q. The Jordan plane and quantum planes with different parameters are non-isomorphic except for quantum planes with inverse parameters, $A_q\simeq A_{q^{-1}}$.
\end{thm}
Theorem~\ref{Jordan} is proved by a direct computation.

\begin{thm}\label{identity} Let $x,y$ be generators of the Jordan plane $A^J$.  For any complex number $p$ and a natural number $k$, the following identity holds:
\be\label{identity formula}
(x+py)^k=(x+(p-k+1)y)(x+(p-k+3)y)(x+(p-k+5)y)\ldots(x+(p+k-1)y).
\ee
\end{thm}

\begin{proof}
For $l\in\N$, we have $xy^l=y^l+ly^{l+1}$. By induction on $k$, we prove that both the right and left hand sides in formula \Ref{identity formula} are equal to
\be  
\sum_{i=0}^k\,{k\choose i}p(p+1)\ldots(p+i-1)y^ix^{k-i}.
\ee
\end{proof}

\subsection{Proof of Theorem~\ref{subbundle}.}\label{proof}

Let $(L_1\T\ldots L_n)^{sing}_l$ be a resonance subspace, otherwise the Theorem is trivial. We prove that for each $m\in (L_1\T\ldots L_n)^{sing}_l$,
\bean\label{commute}
K_i(z)(e(z_1,\dots,z_n))^km=(e(z_1,\dots,z_i+p,\dots,z_n))^kK_i(z)m,
\eean
for $i=1,\dots,n$. 

For $j=1,\dots,n$, we have 
\bean\label{intertwines}
PR_{L_jL_{j+1}}(z_j-z_{j+1})(e(z))^k=(e(z))^kPR_{L_jL_{j+1}}(z_j-z_{j+1}),
\eean
since $e(z)\in Y(gl(2))$ and $PR_{L_jL_{j+1}}(z_j-z_{j+1})$ is an intertwiner of Yangian modules. Therefore, the subspace of conformal blocks is invariant with respect to permutations of variables, $PR_{L_iL_{i+1}}(z_i-z_{i+1})C(z)=C(z_1,\dots,z_{i+1},z_i,\dots,z_n)$.

Consider the case $i=1$. Using \Ref{intertwines}, we get
\be
K_1(z)(e(z))^k=((P_{2,\dots,n,1})^{-1}e(z)P_{2,\dots,n,1})^k K_1(z),
\ee
where $P_{2,\dots,n,1}$ is the permutation operator, $P_{2,\dots,n,1}\,m_1\T\ldots\T m_n=m_2\T\ldots\T m_n\T m_1$, for $m_i\in L_i$, $i=1,\dots,n$.

Since $K_1(z)$ preserves $(L_1\T\ldots\T L_n)_l^{sing}$, it suffices to show that for each $m\in (L_1\T\ldots\T L_n)_l^{sing}$,
\be
((P_{2,\dots,n,1})^{-1}e(z)P_{2,\dots,n,1})^km=(e(z_1+p,z_2,\dots,z_n))^km.
\ee
We have
\be
e(z_1+p,z_2,\dots,z_n)=e(z)+p\,e_\th^{(1)}.
\ee

For $m\in L_1\T\ldots\T L_n$, from \Ref{e(z)}, we get
\bean\label{permuted e(z)}
\lefteqn
{((P_{2,\dots,n,1})^{-1}e(z)P_{2,\dots,n,1})^km=}
\\&&
\left(e(z)+2h_\th^{(1)}e_\th-(2h_\th+(N-2))e_\th^{(1)}+\sum_{j=2}^{N-2}(e_{j1}^{(1)}e_{Nj}-e_{Nj}^{(1)}e_{j1})\right)^k\,m.\notag
\eean

For $s=1,\dots,n$, we have the commutation relations
\bean\label{relations}
e^{(s)}_{Nj}e(z)=e(z)e^{(s)}_{Nj}-e_\th^{(s)}e^{(s)}_{Nj}-\sum_{r=1}^{s-1}2e_\th^{(r)}e^{(s)}_{Nj},
\eean
\be
e_{Nj}e(z)=e(z)e_{Nj}-e_\th e_{Nj},
\ee
where $j=1,\dots,N-1$. Similarly,
\be
e^{(s)}_{j1}e(z)=e(z)e^{(s)}_{j1}-e_\th^{(s)}e_{j1}^{(s)}-\sum_{r=1}^{s-1}2e_\th^{(r)}e^{(s)}_{j1},
\ee
\be
e_{j1}e(z)=e(z)e_{j1}-e_\th e_{j1},
\ee
for $j=2,\dots,N$.

Now, for any $m\in (L_1\T\ldots\T L_n)^{sing}$, we use the commutation relations to move $e_{Nj}, e_{j1}$ and $e_\th$ in formula \Ref{permuted e(z)} to the right and get
\be
((P_{2,\dots,n,1})^{-1}e(z)P_{2,\dots,n,1})^k=(e(z)-(2h_\th+(N-2))e_\th^{(1)})^k\,m,
\ee
since $e_{Nj}m=e_{j1}m=e_\th m=0$ for $j=2,\dots,N-1$.

For $s\in\Z_{\ge 0}$, $m\in (L_1\T\ldots\T L_n)_s$, we have $e(z)m, e_\th^{(1)} m\in (L_1\T\ldots\T L_n)_{s-1}$. Therefore, for each $m\in (L_1\T\ldots\T L_n)_l^{sing}$, the resonance condition \Ref{resonance} implies
\bea
\lefteqn{
(e(z)-(2h_\th+(N-2))e_\th^{(1)})^k\,m=}
\\&&
(e(z)+(p-k+1)e_\th^{(1)})\dots(e(z)+(p+k-3)e_\th^{(1)})(e(z)+(p+k-1)e_\th^{(1)})\,m.
\eea

Notice that the operators $e(z)$ and $e_\th^{(1)}$ define a representation of the Jordan algebra, $e(z)e_\th^{(1)}=e_\th^{(1)}e(z)+e_\th^{(1)} e_\th^{(1)}$, see formula \Ref{relations} for $j=s=1$. 

For $i=1$, formula \Ref{commute} follows from Theorem~\ref{identity}. 

For $i=2,\dots,n$, formula \Ref{commute} follows from formula \Ref{commute} for $i=1$, formula \Ref{intertwines} and the first formula in Lemma~\ref{compatibility}.
$\;\Box$

Notice that $gl(N)$ can be replaced with $sl(N)$. Namely, we consider a tensor product of evaluation modules over $Y(sl(N))$. The same proof shows that the subspace of conformal blocks defined as in Section~\ref{qKZ} is invariant under the qKZ connection.

\end{document}